\documentclass[12pt,twoside]{article}

\usepackage{amsmath,amssymb,amsthm,amscd,color,accents,mathrsfs,url}

\usepackage{stackrel}
\usepackage{enumitem}
\usepackage{euscript}
\usepackage{multirow}
\usepackage{latexsym,amsbsy}
\usepackage[pdftex]{graphicx}
\usepackage{epsfig}
\usepackage{subdepth}
\usepackage[normalem]{ulem}

\usepackage{float}
\usepackage{color}
\usepackage{relsize}
\usepackage[geometry]{ifsym}
\usepackage{pdfpages}
\usepackage[makeroom]{cancel}

\usepackage[showonlyrefs=false]{mathtools}

\usepackage[pagebackref=false,colorlinks=true,linkcolor=blue,citecolor=blue,urlcolor=blue]{hyperref}

\numberwithin{equation}{section}

\theoremstyle{plain}
\newtheorem{theorem}{Theorem}[section]
\newtheorem{proposition}[theorem]{Proposition}

\newtheorem{corollary}[theorem]{Corollary}
\newtheorem{remark}[theorem]{Remark}

\makeatletter
\newcommand*\bigcdot{\mathpalette\bigcdot@{.8}}
\newcommand*\bigcdot@[2]{\mathbin{\vcenter{\hbox{\scalebox{#2}{{\hskip 1pt}$\m@th#1\bullet$}}}}}
\makeatother

\def\C {\mathbb{C}}
\def\E {\mathbb{E}}

\def\R {\mathbb{R}}

\def\vec {\mathop{\mathrm{vec}}}

\def\dd {\hskip2pt\mathrm{d}}

\def\re{\hskip2pt{\rm{Re}}\hskip2pt}

\DeclareMathOperator{\mi}{\mathrm{i}}	

\DeclareMathOperator{\tr}{\mathrm{tr}}
\def\as {\stackrel{\mathrm{a.s.}}{\longrightarrow}}

\newcommand{\cidist}{\stackrel{\mathrm{d}}{\longrightarrow}}

\begin{document}

\title{\vspace{-35pt}
\bf \large
Omnibus Goodness-of-Fit Testing for Distributions \\ on Stiefel Manifolds}

\author{{Dominic Edelmann}\footnote{Division of Biostatistics, German Cancer Research Center, Im Neuenheimer Feld 280, 69120 Heidelberg, Germany.  E-mail address: \href{mailto:dominic.edelmann@dkfz-heidelberg.de}{dominic.edelmann@dkfz-heidelberg.de}.}
\ {and Donald Richards}\footnote{Department of Statistics, Penn State University, University Park, PA 16802, U.S.A.  E-mail address: \href{mailto:richards@stat.psu.edu}{richards@stat.psu.edu}.
\endgraf
\endgraf
\endgraf}
}

\maketitle

\medskip

\date{}

\medskip

\begin{abstract}
In this article, a comprehensive framework for goodness-of-fit testing for distributions on Stiefel manifolds is developed. The approach is based on integrals of the squared differences between empirical and theoretical characteristic functions, yielding test statistics that are consistent against all fixed alternatives. For the Fisher-Bingham family of 
distributions, explicit computable forms of the test statistic are derived. Simplified expressions for important special cases, including the matrix Fisher, matrix Bingham, and uniform distributions are provided. In the case of testing uniformity on hyperspheres, we obtain the complete asymptotic distribution of the test statistic, enabling computationally efficient asymptotic testing. For general Fisher-Bingham distributions, we establish theoretically justified Monte Carlo testing procedures for both simple and composite hypotheses. Simulation studies demonstrate accurate Type I error control and strong power across a wide range of alternatives. The practical relevance of the proposed methodology is illustrated by an application to data on the orbits of comets.
\end{abstract}


\normalsize

\parindent=20pt

\section{Introduction}
\label{sec:introduction}
\setcounter{equation}{0}

In large part, classical statistical procedures were developed for data on Euclidean spaces. By contrast, many modern data sets are inherently non-Euclidean, arising instead on curved geometric spaces where, moreover, standard algebraic operations are inapplicable. Prominent examples of such curved spaces are the Stiefel manifolds $V_{d,p}$, consisting of all $d \times p$ matrices with orthonormal columns.

Data sets for which observations lie in a Stiefel manifold arise in many fields, including medicine \cite{chakraborty2019statistics}, geology \cite{chang1993spherical}, biology \cite{gagliardo2008navigational}, astronomy \cite{lin2017bayesian}, robotics \cite{sola2017quaternion} and other areas \cite{bagyan2024complete,massart2023coordinate,persson2025adaptive}.  Important special cases of Stiefel manifolds are the $(d-1)$-dimensional hyperspheres $S^{d-1}=V_{d,1}$, and the special orthogonal groups $SO(d) \simeq V_{d,d-1}$.

When performing statistical inference or modeling on curved manifolds, goodness-of-fit testing is of the utmost importance. Visual diagnostics are typically infeasible, and deviations from an assumed model can lead to misleading scientific conclusions. Consequently, although there now exists a substantial body of research devoted to goodness-of-fit testing on curved manifolds, most of the extant literature is focused on testing for uniformity on the hypersphere $S^{d-1}$; cf., \cite{garcia2018overview} for an extensive overview. By contrast, the literature on goodness-of-fit testing on general Stiefel manifolds is comparatively sparse \cite{ebner2024unified,ebner2018multivariate,jupp2020measures,jupp2005sobolev,xu2021interpretable}.

In this article, we develop a family of test statistics for goodness-of-fit testing on Stiefel manifolds.  These statistics are derived from integrals that measure the distance between population and empirical characteristic functions. Although this approach is well known in the case of Euclidean data, it seems to have been largely unexplored in the case of curved manifolds, possibly because the calculation of the underlying population characteristic functions were regarded generally as too recondite for practical or theoretical usage. Nevertheless, we establish that these integrals can be analyzed for the purposes of goodness-of-fit testing. 

As a consequence of our analysis, the resulting goodness-of-fit test statistics are proved to be consistent against all alternatives and their asymptotic distributions are shown to be given by an infinite-dimensional Gaussian quadratic form. We also provide explicit formulas for our test statistic in the case of the general Fisher-Bingham distribution on a Stiefel manifold; and as further special cases, we thereby obtain explicit expressions for the goodness-of-fit statistics for the matrix Fisher, matrix Bingham, and the uniform distribution on any Stiefel manifold.

For the case of the uniform distribution on the hypersphere $S^{d-1}$, we derive explicitly the complete asymptotic distribution of the test statistic, obtaining computationally efficient asymptotic test procedures. For the general Fisher-Bingham case, we develop valid resampling-based tests  for both simple and composite hypotheses.

The performance of our test statistics is investigated in a detailed simulation study. These simulations demonstrate that our testing procedures reliably control the nominal Type I error rate and show good power properties over a range of different scenarios for the alternative hypothesis. We provide an application to testing goodness-of-fit for the distribution of a well-known data set on the orbits of comets, and this application also highlights the ability of our method for broader applications to real-world data.

The outline of this article is as follows. In Section~\ref{sec:gof}, we present the general principle underlying our goodness-of-fit approach and provide some important results. Section~\ref{sec:general} states explicit expression for the goodness-of-fit statistics of the general Fisher-Bingham distributions and proposes choices for the hyperparameter matrix that lead to simplification of this test statistic. Section~\ref{sec:special} discusses results for the special cases of the matrix Fisher, matrix Bingham and uniform distributions. In Section~\ref{sec:asy}, we develop an asymptotic test for testing uniformity on the hypersphere. In Section~\ref{sec:mc}, theoretical results are provided for Monte Carlo approaches to testing goodness-of-fit for the general Fisher-Bingham distribution. A detailed simulation study in Section~\ref{sec:sim} and a real-world data example on the orbit of comets in Section~\ref{sec:real} complement our theoretical results. We conclude with a discussion of the article in Section~\ref{sec:discussion}.

\paragraph{Notation.}  For any matrix $M$, we denote by $M'$ the transpose of $M$. For a square matrix $M$, $\tr(M)$ denotes the trace of $M$. For $q \in \mathbb{N}$, $I_q$ denotes the identity matrix in $\R^{q \times q}$; also, $\otimes$ denotes the Kronecker product. We denote by $\vec(\cdot)$ the vec-operator, which is defined for a $d \times p$ rectangular matrix $M = (m_{ij})_{1 \leq i \leq d, 1 \leq j \leq p}$ by  
 $$
 \vec(M) = (m_{11},m_{21},\ldots,m_{d1},m_{12}, \ldots,m_{d2},\ldots,m_{1p},\ldots,m_{dp})',
 $$
the vector obtained by stacking the columns of $M$ from left to right. We also denote by ${\vec}^{-1}_{d \times p}$  the inverse vec-operator which is defined for any $u = (u_1,\ldots, u_{dp})'\in \R^{dp}$ by 
$$
{\vec}^{-1}_{d \times p} (u) = 
\begin{bmatrix}
u_{1} &u_{d+1}&\cdots&u_{d(p-1)  +1}\\
u_{2}&u_{d+2}&\ddots&u_{d(p-1)   +2}\\
\vdots&\ddots & \ddots &\vdots\\
u_{d}&u_{d+p}&\cdots&u_{dp}
\end{bmatrix} 
\in\mathbb{R}^{d \times p},
$$
which is the $d \times p$ matrix that is obtained by filling up the columns of a $d \times p$ matrix from left to right. 
Further, we use the notation $\as$ to denote almost sure convergence, and convergence in distribution will be denoted by $\stackrel{\mathrm{d}}{\longrightarrow}$. 

Let \( \Omega \subseteq \mathbb{R}^{d \times p} \), and let \( \mathcal{F} \) denote the Borel \( \sigma \)-algebra on \( \Omega \). Given a measure \( \mu \) on \( (\Omega, \mathcal{F}) \) and $p \ge 1$, we denote by $L^p(\Omega,\mu)$ the space of (equivalence classes of) measurable functions $f: \Omega\to \C$ such that $\int_{\Omega} |f|^p \dd \mu(x) < \infty$.

\section{Goodness-of-fit testing based on the empirical \\ characteristic function} \label{sec:gof}


Consider a random sample $X_1,\ldots,X_n \in  \R^{d \times p}$ , drawn from a population with distribution $P$. We wish to test the null hypothesis $H_0: P = P_0$, for some given distribution $P_0$ on $\R^{d \times p}$, against the alternative hypothesis $H_1: P \neq P_0$. 



Let $\mi = \sqrt{-1}$, and denote the population characteristic function under $H_0$ by
$$
\varphi_0(t) := \E \exp(\mi {\tr(t'X)}) 
= \int_{{\R^{d \times p}}} \exp(\mi \tr (t'x)) \dd P_0 (x), \ 
$$
$t \in \R^{d \times p}$.
Further, the empirical characteristic function of the sample $X_1,\ldots,X_n$ is defined as 
$$
\varphi_n(t) = \frac{1}{n} \sum_{j=1}^n \exp(\mi {\tr(t'X_j)}) .
$$

Let $w: \R^{d \times p} \to [0,\infty)$ be a nonnegative weight function such that $w \in L^1(\R^{d \times p})$. Using this weight function, we define the goodness-of-fit test statistic
\begin{equation} \label{eq:teststat}
D_n = \int_{\R^{d \times p}} |\varphi_n(t) - \varphi_0(t)|^2 w(t) \dd t.
\end{equation}
Without loss of generality, we normalize the weight function to have $L^1$-norm $1$, so we assume throughout that
\begin{itemize}
    \item[(C.1)] $w: \R^{d \times p} \to \R$ is a Lebesgue-measurable probability density function.
\end{itemize}
By applying to Eq.~\eqref{eq:teststat} the expansion 
\begin{align} \label{eq:Dn_expansion}
|\varphi_n(t) - \varphi_0(t)|^2 &= \big(\varphi_n(t) - \varphi_0(t)\big) \big(\overline{\varphi_n(t)} - \overline{\varphi_0(t)}\big) \nonumber \\
&= |\varphi_n(t)|^2 - 2 \re \varphi_n(t) \overline{\varphi_0(t)} + |\varphi_0(t)|^2,
\end{align}
and noting that the norm of any characteristic function is less than or equal to $1$, then it follows that $0 \le D_n \le 4 $.  Further, 
we obtain from Eq.~\eqref{eq:teststat} a decomposition of $D_n$ into three terms,  Further, 
\begin{align} \label{eq:Dn_Un_expansion}
D_n 
&= U_{n,1} - 2 U_{n,2} + U_{3},
\end{align}
where 
\begin{equation}
\label{eq:Un_123}
\begin{aligned}
U_{n,1} &:= 
\int_{\R^{d \times p}} |\varphi_n(t)|^2 w(t) \dd t, \\
U_{n,2} &:= 
\re \int_{\R^{d \times p}} \varphi_0(t) \overline{\varphi_n(t)} w(t) \dd t, \\
U_3 &:= 
\int_{\R^{d \times p}} |\varphi_0(t)|^2 w(t) \dd t.
\end{aligned}
\end{equation}
The statistic $U_{n,1}$ depends on the sample $X_1,\ldots,X_n$ only, whereas the statistic $U_{n,2}$ depends on both the random sample and the density function $f_0$. The term $U_{3}$ depends on $f_0$ only, which explains why we omit the index $n$ for that term. 

Define the Fourier transform of $w(t)$, viz., 
$$
\widehat{w}(x) = \int_{\R^{d \times p}} \exp(\mi \tr (t'x)) w(t) \dd t,
$$
$x \in \R^{d \times p}$.  Note that, since $w$ is a probability density function, then $\widehat{w}$ is positive definite.  Moreover, $\widehat{w}$ is Hermitian, i.e., $\widehat{w}(x) = \overline{\widehat{w}(-x)}$, $x \in \R^{d \times p}$.  

Also define \begin{equation} \label{eq:capitalW}
\widehat{W}_0(x) = \int_{\R^{d \times p}} \widehat{w}(x-u) \dd P_0 (u),
\end{equation}
$x \in \R^{d \times p}$.  Then, $\widehat{W}_0$ is a generalized convolution. In particular, if $P_0$ possesses a density function $f_0$ with respect to the Lebesgue measure on $\R^{d \times p}$ then Eq.~\eqref{eq:capitalW} simplifies to
$$
\widehat{W}_0(x) = (\widehat{w} * f_0)(x).
$$

We now show that each of the terms $U_{n,1}$, $U_{n,2}$, and $U_{3}$ 
can be expressed in terms of $\widehat{w}$, the Fourier transform of $w$.  

\begin{theorem} \label{th:uns}
Assume condition (C.1). Then,
\begin{align}
 U_{n,1} &=  \frac{1}{n^2} \sum_{j,k=1}^n \widehat{w}(X_j-X_k), \label{eq:un1} \\
 U_{n,2} &= \frac{1}{n} \re \sum_{j=1}^n \widehat{W}_0(X_j) , \label{eq:un2} \\
 \intertext{and}
 U_{3} &= \int_{\R^{d \times p}} \widehat{W}_0(v) \dd P_0(v) = {{\int_{\R^{d \times p}}}} \int_{\R^{d \times p}} \widehat{w}(u-v) \dd P_0(u) \dd P_0(v). \label{eq:un3}
\end{align}
\end{theorem}

\medskip

Theorem \ref{th:uns} implies that $D_n$ is a V-statistic with kernel function 
$$
 h(x,y) =  \int_{\R^{d \times p}} \widehat{W}_0(v) \dd P_0(v) + \re \big(\widehat{w}(x-y) - \widehat{W}_0(x) - \widehat{W}_0(y)\big).
$$
As shown within the proof of Proposition \ref{prop:clt} in the Supplementary Material, the kernel $\widehat{w}(x-y)$ is positive definite.  Therefore $h(x,y)$ also is positive definite, and $D_n$ is degenerate of order $1$ if $P=P_0$. Consequently, the asymptotic distribution of $D_n$ can now be derived using results from the theory of {U}-statistics \cite{boroskikh2020u}.




\begin{theorem} \label{prop:clt}
Assume condition (C.1), and consider a random sample $X_1,\ldots,X_n$ $\in \R^{d \times p}$, drawn from a population with distribution $P_0$. 
Then 
$$
n D_n \stackrel{\mathrm{d}}{\longrightarrow} \sum_{i=1}^\infty \mu_i Q_i^2 
$$
as $n \to \infty$, where $Q_1^2,Q_2^2,\ldots$ are mutually independent, identically distributed chi-squared random variables with one degree-of-freedom; $\mu_1 \geq \mu_2 \geq \cdots$, and $\mu_i \geq 0$ for all $i=1,2,3,\ldots$; $\sum_{i=1}^\infty \mu_i < \infty$; and the $\mu_i$ are the eigenvalues of the trace-class integral operator $H: L^2(\R^{d \times p}, P_0) \to L^2(\R^{d \times p}, P_0)$ such that, for $f \in L^2(\R^{d \times p}, P_0)$, 
$$
 (H f)(x) = \int_{\R^{d \times p}} h(x,y) f(y) \dd P_0 (y), \qquad x \in \R^{d \times p}.
$$
\end{theorem}

\medskip

The strong consistency of $D_n$ can be derived using the Strong Law of Large Numbers for V-statistics \cite{gine1992marcinkiewicz}.

\begin{proposition}
Assume condition (C.1) and consider a random sample $X_1,\ldots,X_n \in \R^{d \times p}$, drawn from a population with distribution $P$.
 \begin{enumerate}
  \item[(i)] If $P = P_0$ then $n D_n \as 0$ as $n \to \infty$.  
  \item[(ii)] If $w$ is positive almost everywhere and $P \neq P_0$ then ${{D_n \as c}}$ as $n \to \infty$, where $0 < c \le 4$.
 \end{enumerate}
\end{proposition}

We remark that, as alternatives to deriving the test statistic $D_n$ using a characteristic function-based approach, we can also derive equivalent test statistics using the maximum mean discrepancy  \cite{gretton2012kernel} or the generalized energy distance \cite{edelmann2022regression,rizzo2016energy,sejdinovic2013equivalence}.

\section{The general Fisher-Bingham distribution} \label{sec:general}

In this section, we consider the problem of goodness-of-fit testing for the Fisher-Bingham distribution on the Stiefel manifold
    $$
    V_{d,p} := \{x \in \R^{d \times p} : x' x  = I_p\},
    $$ 
$d \geq p$.  
For $A \in \R^{d \times p}$ and a symmetric matrix $B \in \R^{dp \times dp}$, the density of the Fisher-Bingham distribution on $V_{d,p}$ with respect to the normalized Haar measure can be written as 
\begin{equation} \label{eq:gendensity2}
  g_0 (x; A, B) = [\Psi(A,B)]^{-1} \exp \left( \tr(A' x) + \vec(x)' \, B \, \vec(x) \right),
\end{equation}
$x \in V_{d,p}$, where the normalizing constant $\Psi (A,B)$ is given by
\begin{equation} \label{eq:normconstant}
    \Psi (A, B) = \int_{V_{d,p}} \exp \left(\tr(A' v)  + \vec(v)' \, B \, \vec(v) \right) \dd \eta(v), 
\end{equation}
and $\dd \eta$ is the normalized Haar measure on $V_{d,p}$.  
We will also use the alternative representation, 
\begin{equation} \label{eq:gendensity}
g_0 (x;{\vec}^{-1}_{d \times p} (a),B) = [\Psi_{\vec}(a, B)]^{-1} \exp \left(a' \vec(x) + \vec(x)' \, B \, \vec(x) \right) ,
\end{equation}
with $ \Psi_{\vec} (\vec(A), B) = \Psi(A,B)$.

The normalizing constant $\Psi (A, B)$ is, in general, not available in closed form, cf., \cite{kume2013saddlepoint}; however, several important special cases admit explicit expressions. One notable example is the matrix Langevin distribution that arises as a special case of \eqref{eq:gendensity2} with $B = 0$, for which the normalizing constant is \cite[p. 31]{chikuse2012statistics},
 \begin{equation} \label{eq:normconstanthyp}
    \Psi(A,0) = {}_0F_1(\tfrac12 d; \tfrac14 A' A),
\end{equation}
where ${}_0F_1$ is the generalized hypergeometric function of matrix argument \cite{gross1987special,muirhead1982aspects}.

For the problem of testing that a random sample on $V_{d,p}$ is drawn from a Fisher-Bingham distribution \eqref{eq:gendensity2}, we apply the statistic $D_n$ arising from Eq.~\eqref{eq:teststat} with the weight function
\begin{equation} \label{eq:weight}
    w(t;\Omega) = (2 \pi)^{- dp/2} \det(\Omega)^{- 1/2} \exp \left(- \tfrac{1}{2} \vec(t)'  \, \Omega^{-1} \, \vec(t) \right),
\end{equation}
 where $\Omega$ is a positive definite (symmetric) matrix in $\R^{dp \times dp}$. The Fourier transform of $w$ is well known to be
\begin{equation}
    \widehat{w}(x; \Omega) = \exp \left(- \tfrac{1}{2} \vec(x)'  \, \Omega \, \vec(x) \right). \label{eq:what}
\end{equation}
By applying Theorem \ref{th:uns}, we derive the following representations for $U_{n,1}$ and $U_{n,2}$, and we emphasize the dependence of each term on $\Omega$ by using the notation $U_{n,1; \Omega}$ and $U_{n,2; \Omega}$, respectively.

\begin{theorem} \label{th:general}
Let $X_1,\ldots,X_n \in  \R^{d \times p}$ denote a random sample on ${V}_{d,p}$. Set $f_0(x) = g_0 (x; A,B)$ and choose the weight function $w$ as in Eq.~\eqref{eq:weight}. Then the first two terms in Theorem \ref{th:uns} can be expressed as
\begin{align*}
    U_{n,1; \Omega} &= \frac{1}{n^2} \sum_{j,k=1}^n \exp \left(- \tfrac{1}{2} \vec(X_j - X_k)'  \, \Omega \, \vec(X_j-X_k) \right), \\
    U_{n,2; \Omega} &= \frac{1}{n \, \Psi (A, B)} \sum_{j=1}^n \exp \left(- \tfrac{1}{2} \vec(X_j)'  \, \Omega \, \vec(X_j) \right) \, \Psi_{\vec}(a+ \Omega \vec(X_j), B-\tfrac{1}{2} \Omega).
\end{align*}
\end{theorem}

In Section \ref{subsec:uniform}, we calculate the constant term $U_{3}$ explicitly for the special case of the uniform distribution on the Stiefel manifold. Although it can be difficult in general to obtain a simple expression for $U_{3}$, we will show in Section~\ref{sec:mc} that $U_{3}$ can be evaluated using Monte Carlo methods or other numerical procedures. Alternatively, we can circumvent the calculation of $U_3$ by constructing Monte Carlo-based tests using the equivalent test statistic,
\begin{equation} \label{eq:equivalent}
   \widetilde{D}_n =  U_{n,1; \Omega} - 2 \, U_{n,2; \Omega}.
\end{equation}

As for the statistic $U_{n,2; \Omega}$, its only non-trivial terms are normalizing constants of the form $\Psi_{\vec}(s,R)$ for some $s \in \R^{d p}$ and $R \in \R^{dp \times dp}$. For general $s$ and $R$, no closed-form expressions are available, however several numerical approaches for computing $\Psi_{\vec}(s,R)$ have been developed \cite{chen2021maximum,kume2013saddlepoint,kume2018exact}. 

Depending on the desired precision for these normalizing constants, it may be computationally expensive to perform the test, particularly when a Monte Carlo approach is used.  However, this computational burden can be alleviated by choosing a matrix $\Omega$ that leads to a simplification of the terms in Theorem \ref{th:general}. In particular, let $\Lambda$ be a symmetric $p \times p$ matrix such that 
\begin{equation} \label{eq:omegalambda}
\Omega =  2 \, (B + \Lambda \otimes I_d) 
\end{equation}
is positive definite. Such a matrix $\Lambda$ can be found through direct calculation.  Alternatively, denote by $\lambda_{\min}(S)$ the smallest eigenvalue of a symmetric matrix $S$; then by Weyl's inequality \cite[p.~239]{horn2012matrix}, 
\begin{align*}
\lambda_{\min}(\Omega) &= 2 \lambda_{\min}(B + \Lambda \otimes I_d) \\
&\ge 2 \, [\lambda_{\min}(B) + \lambda_{\min}(\Lambda \otimes I_d)] 
= 2 \, [\lambda_{\min}(B) + \lambda_{\min}(\Lambda)],
\end{align*}
hence $\Omega$ in Eq.~\eqref{eq:omegalambda} is positive definite whenever $\lambda_{\min}(B) + \lambda_{\min}(\Lambda) > 0$.  

Now setting $\Omega = 2 \, (B + \Lambda \otimes I_d) $, the terms $\Psi_{\vec}(a+ \Omega \vec(X_j), B-\tfrac{1}{2} \Omega)$ reduce to constant multiples of $\Psi_{\vec}(a+ \Omega \vec(X_j), 0)$, which can be evaluated in terms of generalized hypergeometric functions; cf. Eq.~\eqref{eq:normconstant}.

\begin{proposition} \label{prop:stat1}
Let $\Lambda$ be a symmetric $p \times p$ matrix such that $\Omega = 2(B + \Lambda \otimes I_d)$ is positive definite. Then the terms in Theorem \ref{th:general} can be written as
\begin{align*}
U_{n,1; 2 (B+\Lambda \otimes I_d)} &= \frac{\exp(- 2 \tr(\Lambda))}{n^2} \\
& \quad \times \sum_{j,k=1}^n \left( \exp \left(- \vec(X_j - X_k)'  \, B \, \vec(X_j-X_k) \right) \exp(2 \tr (\Lambda X_j' X_k)) \right), \\
\intertext{and}
U_{n,2; 2 (B+\Lambda \otimes I_d)} &=  \frac{\exp(- 2 \tr(\Lambda))  }{n } [\Psi(A, B)]^{-1} \\
& \qquad\qquad\qquad \times \sum_{j=1}^n \exp (- \vec(X_j)'  \, B \, \vec(X_j) ) \, {}_0F_1(\tfrac12 d; \tfrac14 {{Z}}_j' {{Z}}_j),
\end{align*}
where ${Z}_j = A + 2 X_j \Lambda +2 {\vec}^{-1}_{d \times p} (B \, \vec(X_j))$, $j=1,\ldots,n$.
\end{proposition}

Due to the availability of efficient algorithms for evaluating the ${}_0F_1$ functions \cite{koev2006efficient}, every summand in the expression for $U_{n,2 ;2 (B+\Lambda \otimes I_d)}$ can be computed rapidly. This allows us to use very precise algorithms for calculating the single normalizing constant $\Psi_{\vec}(A,B)$ without substantial computational effort.

Another interesting choice for $\Omega$ is $\Omega = 2 (\Lambda \otimes I_d)$, where $\Lambda$ is a positive definite $p \times p$ matrix. As shown in the following result, this choice leads to compact expressions for $U_{n,1}$ and $U_{n,2}$.

\begin{proposition} \label{prop:stat2}
Let $\Omega = 2 \Lambda \otimes I_d$, where $\Lambda$ is a positive definite ${p \times p}$ matrix. Then the terms in Theorem \ref{th:general} can be written as
\begin{align*}
U_{n,1; 2 \, \Lambda \otimes I_d} & = \frac{\exp(- 2 \tr(\Lambda))}{n^2} \sum_{j,k=1}^n\exp(2 \tr (\Lambda X_j' X_k)), \\
\intertext{and}
U_{n,2; 2 \, \Lambda \otimes I_d} &=  \frac{\exp(- 2 \tr(\Lambda))  }{n } [\Psi (A, B)]^{-1} \,\sum_{j=1}^n \Psi(A+2 \, X_j \Lambda, B).
\end{align*}
\end{proposition}


\section{Special cases} \label{sec:special}

\subsection{The matrix Fisher distribution}\label{subsec:matrixfisher}

The matrix Fisher distribution - which is also referred to as the matrix Langevin or matrix von Mises-Fisher distribution - is obtained by setting $B=0$ in the density function of the Fisher-Bingham distribution in Eq.~\eqref{eq:gendensity}. Taking into account the considerations of the previous section, a natural choice for $\Omega$ is 
\begin{equation} \label{eq:langevinomega}
\Omega = 2 \Lambda \otimes I_d,
\end{equation}
where $\Lambda$ is a positive definite $p \times p$ matrix.  With this specification of $\Omega$, the terms $U_{n,1}$ and $U_{n,2}$ admit closed-form expressions involving only elementary terms and generalized hypergeometric functions.

\begin{corollary} \label{cor:fisher}
Consider the matrix Fisher distribution, and let $\Omega$ be defined as in \eqref{eq:langevinomega}. Then the terms in Theorem \ref{th:general} can be written as
\begin{align*}
U_{n,1 ; 2 \Lambda \otimes I_d } &= \frac{\exp(- 2 \tr(\Lambda))}{n^2} \sum_{j,k=1}^n \exp(2 \tr(\Lambda X_j' X_k)),  \\
\intertext{and}
U_{n,2 ; 2 \Lambda \otimes I_d}  &= \frac{\exp(- 2 \tr(\Lambda))  }{n } [ {}_0F_1(\tfrac12 d; \tfrac14 A' A)]^{-1} \sum_{j=1}^n {}_0F_1 \left(\tfrac12 d; \tfrac14 Z_j' Z_j \right),
\end{align*}
where $Z_j = A + 2 X_j \Lambda$, $j=1,\ldots,n$.  
\end{corollary}

\subsection{The matrix Bingham distribution}\label{subsec:matrixbingham}

With symmetric matrices $D \in \R^{p \times p}$ and $E \in \R^{d \times d}$, the density of the matrix Bingham distribution on $V_{d,p}$ with respect to the normalized Haar measure is 
\begin{equation}
\big[{}_{0}F_0(D,E)\big]^{-1} \exp \left(\tr(D x' E x) \right) ,
\end{equation}
$x \in V_{d,p}$, where the function ${}_{0}F_0(D,E)$ is a generalized hypergeometric function of two matrix arguments \cite[p. 108]{chikuse2012statistics}.  By \cite[p.~76, Lemma 2.2.3]{muirhead1982aspects}
and the fact that $E$ is symmetric, 
$$
\tr(D x' E x)  = \vec(x)' (D \otimes E') \vec(x) = \vec(x)' (D \otimes E) \vec(x),
$$
Hence the matrix Bingham distribution can be written as a general Fisher-Bingham distribution with density $g_0 (x;0,D \otimes E)$.  

As in Proposition \ref{prop:stat1}, we choose
$$
 \Omega = \Omega(D,E) = 2 \, (D \otimes E + \Lambda \otimes I_d),
$$
where $\Lambda$ is a symmetric matrix that defines a positive definite $\Omega$. This yields again closed-form expressions for $U_{n,1}$ and $U_{n,2}$  involving only elementary terms and generalized hypergeometric functions.

\begin{corollary}
Consider the matrix Bingham distribution, i.e., set $A=0$ and $B = D \otimes E$. Let $\Lambda$ be a symmetric matrix  in $\mathbb{R}^{p \times p}$, such that $\Omega = 2(B + \Lambda \otimes I_d)$ is positive definite. Then the terms in Theorem \ref{th:general} can be written as
\begin{align*}
U_{n,1; \Omega (B, \Lambda) } &= \frac{\exp(- 2 \tr(\Lambda))}{n^2} \\ 
&\quad\qquad \times \sum_{j,k=1}^n \exp \left(- \tr(D (X_j-X_k)' E (X_j-X_k) )\right) \exp(2 \tr (\Lambda X_j' X_k)),
\end{align*}
and 
$$
U_{n,2; \Omega (B, \Lambda) } =  \frac{\exp(- 2 \tr(\Lambda))  }{n } {}_0 F_0 (D, E)^{-1} \,\sum_{j=1}^n \exp \left(-\tr(D X_j' E X_j) \right) \, {}_0F_1(\tfrac12 d; \tfrac14 Z_j' Z_j),
$$
where $Z_j = 2 (X_j \Lambda + E X_j D)$, $j=1,\ldots,n$.
\end{corollary}

\subsection{The uniform distribution} \label{subsec:uniform}

The uniform distribution is, by far, the most studied distribution on the Stiefel manifold \cite{mardia1977uniform, chen2021maximum,iwashita2017test}. The density of the uniform distribution is naturally given by the normalized Haar measure on $V_{d,p}$ and it is obtained by setting $A=0$ and $B=0$ in Eq.~\eqref{eq:gendensity}. Since the uniform distribution is a special case of both the matrix Langevin and the matrix Bingham distribution, we can obtain representations of $U_{n,1}$ and $U_{n,2}$ in simple closed-form expressions. Moreover, it follows by a straightforward symmetry consideration that $U_{3} = U_{n,2}$. 

\begin{corollary} \label{cor:testing_unif}
Consider the uniform distribution on the Stiefel manifold $V_{d,p}$, i.e., set $A=0$ and $B=0$. Let $\Omega = 2 \Lambda \otimes I_d$, where $\Lambda$ is a positive definite ${p \times p}$ matrix. Then the terms in Theorem \ref{th:general} can be written as
\begin{align*}
 U_{n,1 ; 2 \Lambda \otimes I_d} &= \frac{\exp(- 2 \tr(\Lambda))}{n^2}  \, \sum_{j,k=1}^n  \exp(2 \tr (\Lambda X_j' X_k)),  \\
\intertext{and}
 U_{n,2 ; 2 \Lambda \otimes I_d} &=  \exp(- 2 \tr(\Lambda)) \,  {}_0F_1(\tfrac12 d;  \Lambda ^{2}), 
\end{align*}
Moreover, the term $U_{3}$ in Theorem \ref{th:uns} can be evaluated in closed form, 
$$
U_{3 ; 2 \Lambda \otimes I_d} = U_{n,2 ; 2 \Lambda \otimes I_d} = \exp(- 2 \tr(\Lambda)) \, {}_0F_1(\tfrac12 d;  \Lambda^{2}),
$$
and we also have 
\begin{equation} \label{eq:Dn}
D_n = \frac{\exp(- 2 \tr(\Lambda))}{n^2} \sum_{j,k=1}^n \left[ \, \exp(2 \tr(\Lambda X_j' X_k)) - {}_0F_1(\tfrac12 d;  \Lambda^2)\right].
\end{equation}
\end{corollary}

\medskip

We will also consider the equivalent test statistic
\begin{equation} \label{eq:Sn}
E_n = \exp(2 \tr(\Lambda)) D_n = \frac{1}{n^2} \sum_{j,k=1}^n \, \exp(2 \tr (\Lambda X_j' X_k)) - {}_0F_1(\tfrac12 d; \Lambda^2).
\end{equation}

\medskip

In testing for uniformity on the Stiefel manifold $V_{d,p}$ it is interesting to ascertain the behavior of the statistic $E_n$ or $D_n$, for extreme values of $\Lambda$, with $n$, $d$, and $p$ held fixed. 

\begin{remark}
\label{rem:limiting_test_statistic}
    Suppose, in Corollary \ref{cor:testing_unif}, that $\Lambda = \lambda I_p$. Then,
    $$
   \lim_{\lambda \to 0} \lambda^{-1} E_n = 2 \tr(\bar{X}'\bar{X}),
    $$
    which is a constant multiple of the generalized Rayleigh test statistic.
\end{remark}

\begin{remark}
Suppose, in Corollary \ref{cor:testing_unif}, that the distribution of $X_1$ is absolutely continuous with respect to the Haar measure on $V_{d,p}$. Then, denoting by $\tau_{\min}(\Lambda)$ the smallest eigenvalue of $\Lambda$, we obtain
\begin{equation} \label{eq:Dn_limit}
\lim_{\tau_{\min}(\Lambda) \to \infty} \exp(-2 \tr (\Lambda)) E_n = \lim_{\tau_{\min}(\Lambda) \to \infty} D_n = \frac{1}{n},
\end{equation}
i.e., $D_n$ becomes degenerate.
\end{remark}


\section{The asymptotic distribution of the test statistic} \label{sec:asy}

In this section we provide the asymptotic distribution of the test statistic $E_n$, as $n \to \infty$, for the case in which $P_0$ is the uniform distribution on the hypersphere $S^{d-1}$.  In this setting, Eq.~\eqref{eq:Sn} reduces to
\begin{equation} \label{eq:Sn_hypersphere}
E_n = \frac{1}{n^2} \sum_{j,k=1}^n \exp(2 \lambda \, X_j' \, X_k) - {}_0F_1(\tfrac12 d; \lambda^{2}),
\end{equation}
where $\lambda > 0$ and ${}_0F_1(\tfrac12 d;\lambda^2)$ is a classical generalized hypergeometric function with scalar argument $\lambda^2$.  

In order to ascertain the asymptotic distribution we will require the modified Bessel function of the first kind of order $\nu$ (\href{https://dlmf.nist.gov/10.25.E2}{Olver, \textit{et al.} \cite{NIST:DLMF}}), denoted by $I_{\nu}(\cdot)$.

\begin{theorem} \label{th:asy}
Let {$X_1,\ldots,X_n$} be a random sample from the uniform distribution on the hypersphere $S^{d-1}$. Then, as $n \to \infty$,
$$
n E_n \cidist \sum_{k=1}^\infty \mu_k {{Q}}_k^2,
$$
where {{$Q_1^2,Q_2^2,\ldots$}} are independently distributed; for each $k=1,2,3,\ldots$, $Q_k^2$ is chi-squared distributed with 
$$
d_k =\binom{d+k-1}{k} - \binom{d+k-3}{k-2}
$$ 
degrees-of-freedom; and 
$$
\mu_k = \Gamma(d/2) \, \lambda^{-(d-2)/2} \, I_{(d/2)+k-1}(2 \lambda) > 0. 
$$
Further, the sequence $\{\mu_k: k \geq 1\}$ is strictly decreasing.
\end{theorem}

We note that an equivalent statistic and method for testing uniformity on the hypersphere was derived in \cite{fernandez2023new}, and it was also shown there that this test belongs to the class of Sobolev tests developed in \cite{gine1975invariant}.

As a special case of Remark \ref{rem:limiting_test_statistic}, it follows that, for fixed $n$, $d$, and $p$, the statistic $\lambda^{-1} E_n$ converges to a constant multiple of the classical Rayleigh test statistic as $\lambda \to 0$.  Also, $D_n \to n^{-1}$ as $\lambda \to \infty$.

\section{Monte Carlo tests} \label{sec:mc}

\subsection{Testing simple hypotheses} \label{subsec:simple}

\paragraph{A sampling-based test.} Consider a random sample $X_1, \ldots,X_n$ from $X$, a random variable on the Stiefel manifold $V_{d,p}$.  We first consider the problem of testing simple null hypotheses of the form
\begin{align*}
H_0: \ &\textrm{The random variable } X \textrm{ follows a Fisher-Bingham distribution $g(\cdot \,;A,B)$ with} \\
\ &\textrm{given parameters $A$ and $B$}
\end{align*}
against the alternative $H_1$:~$H_0 \textrm{ is not valid}$.

A direct approach for testing in this setting is to generate $K$ mutually independent random samples, each of size $n$, from $g(\cdot \,;A,B)$ and calculate the test statistic in Eq.~\eqref{eq:equivalent} for each $k \in \{1,\ldots,K\}$, which we denote by $\widetilde{D}^*_{kn}$. Then, by construction, for all $u \in \mathbb{R}$ 
\begin{equation} \label{eq:pvalbs}
\xi_{K,u} = \frac{1}{K+1} \left(1+\sum_{k=1}^K 1 (\widetilde{D}^*_{kn} \geq u)\right)
\end{equation}
satisfies, for $K \to \infty$,
$$
\xi_{K,u} \as P_0(D_n \geq u).
$$

Despite its simplicity, there are two difficulties with this testing procedure. First, although methods for sampling from the Fisher-Bingham distribution on the Stiefel manifold $V_{d,p}$ have been derived \cite{hoff2009simulation}, these algorithms can be slow, particularly when $p$ and $d$ are large. Second, the evaluation of each of the $K+1$ test statistics \textit{via} the form in Theorem \ref{th:general} requires the calculation of $(K+1)\,n +1$ normalizing constants of the Fisher-Bingham distribution (however, we note that this computational burden is substantially reduced in special cases such as the matrix Fisher distribution or when $\Omega$ is suitably chosen, cf. Proposition \ref{prop:stat1}).  

\paragraph{A bootstrap-based test.} Another way to construct a Monte Carlo-based testing procedure is to follow the bootstrap approach in \cite{arcones1992bootstrap}. For this approach, we draw $K$ samples of size $n$ with replacement from $X_1,\ldots,X_n$. Denote the $k$-th bootstrap sample by $X_{k1}^*,\ldots,X_{kn}^*$. Then, for each $k \in \{1,\ldots,K\}$ calculate the statistic
\begin{align}
{D}_{kn}^* = \frac{1}{n^2} \Bigg[ \sum_{i,j=1}^n h(X_{ki}^*, X_{kj}^*) - \frac{1}{n} & \sum_{l = 1}^n h(X_{ki}^*, X_{l}) \nonumber \\ &- \frac{1}{n} \sum_{l = 1}^n h(X_{kj}^*, X_{l}) + \frac{1}{n^2} \sum_{l,m = 1}^n h(X_{l}, X_{m}) \Bigg].
\end{align}

If $H_0$ is valid then, by \cite[Theorem 3.5]{arcones1992bootstrap}, for each $u \in \mathbb{R}$,
\begin{align} \label{eq:bs1conv}
\big|P(n \, {D}_{kn}^* > u|X_1,\ldots,X_n) - P(n \, D_n \geq u)\big| \longrightarrow 0,
\end{align}
as $n \to \infty$. Hence an asymptotically valid test for $H_0$ can be established analogous to the sampling-based test, cf. Eq.~\eqref{eq:pvalbs}.

A disadvantage of the bootstrap test is that it requires the calculation of the term $U_{3}$, which is usually not available in closed form. Thus in our simulations we calculate $U_{3}$ \textit{via} a Monte Carlo approach. Observing, that, under $H_0$, as $N \to \infty$,
\begin{equation} \label{eq:u3conv}
-D_{N}^* \as U_3,
\end{equation}
the constant $U_{3}$ can be consistently estimated by drawing a large number of samples $X_1^*,\ldots,X_N^*$ from $P_0$ and evaluating $-D_{N}^*$ for this sample. Since the convergence rate in \eqref{eq:u3conv} is $N^{-1}$, we then require that $N/n \to \infty$ as $n \to \infty$ in order for \eqref{eq:bs1conv} to hold when using this Monte Carlo approach.

\subsection{Testing composite hypotheses} \label{subsec:comp}

A more complex problem arises in testing composite hypotheses such as,
$$
H_0: X \text{ follows a matrix Fisher distribution}
$$
or, in explicit form,
\begin{equation} \label{eq:comp}
H_0: P^X  \in \{P_{g(A,B)} | A \in \R^{d \times p}, \, B = 0 \}
\end{equation}
where $P^X$ denotes the distribution of $X$ and $P_{g(A,B)}$ denotes the distribution corresponding to the density function $g(\cdot;A,B)$. Bootstrap procedures for testing composite hypotheses can be derived along the lines of \cite[Section 3]{leucht2009consistency}. 

Consider the hypothesis \eqref{eq:comp}. Given a random sample $X_1,\ldots,X_n$, we proceed as follows:
\begin{enumerate}
\item Assume that $X_1,\ldots,X_n$ follow a matrix Fisher distribution and compute an estimator $\widehat{A}_n$ of $A$ (e.g., \textit{via} \cite[Theorem 2]{jupp1979maximum}).
\item Calculate the test statistic $D_n$ based on setting $A = \widehat{A}_n$ in Corollary \ref{cor:fisher}. 
\item Sample $K$ i.i.d. samples of size $n$ from $P_{g(\widehat{A},0)}$.
\item For the $k$th random sample, calculate the corresponding estimator ${A}^*_{kn}$, $1 \le k \le K$.
\item For the $k$th random sample, calculate $D_{kn}^*$, the test statistic based on setting $A={A}^*_{kn}$ in Corollary \ref{cor:fisher}, $1 \le k \le K$. 
\end{enumerate}

\noindent
Under certain assumptions, this procedure yields an asymptotically valid test for the hypothesis in \eqref{eq:comp}. 
We now show, employing arguments similar to \cite[Section 3]{leucht2009consistency}, that the computational bootstrap is consistent for a broad class of estimators. The key requirement is that the estimator admits an asymptotic linear representation — a condition satisfied by most standard estimators, including MLEs and M-estimators.

\begin{theorem} \label{th:compbs}
Consider an i.i.d.\ sample $X_1,\ldots,X_n$ drawn from $P_{g(A,0)}$, and
write $a = \vec(A)$. Let $\widehat{a}_n = \vec(\widehat{A}_n)$ be an
estimator of $a$, and assume that, as $n \to \infty$,
$$
\widehat{a}_n \as a.
$$
Assume further that $\widehat{a}_n$ admits the asymptotic expansion
$$
\widehat{a}_n - a = \frac{1}{n} \sum_{i=1}^n l(a, X_i) + o_P(n^{-1/2}),
$$
where $l : \R^{dp} \times \R^{dp} \to \R^{dp}$ is a continuous function
satisfying $E[l(b, X_i)] = 0$ for every $b \in \R^{dp}$, together with
$\sup_{b \in N} E[\|l(b, X_i)\|^2] = C_N < \infty$ for some
neighborhood $N$ of $a$. Then, as $n \to \infty$,
$$
\sup_{u \in \mathbb{R}} \big| P(n\, D_{kn}^* > u \mid X_1, \ldots, X_n)
- P(n\, D_n \geq u) \big| \as 0.
$$
\end{theorem}

As for the standard bootstrap approach, the $U_3$ terms are typically not available in closed form, and then are computed using, e.g., the Monte Carlo approach described before. Using the same arguments as before, it can be shown that Theorem \ref{th:compbs} remains valid using this Monte Carlo approach under the assumption that $N/n \to \infty$ as $n \to \infty$.

\section{Simulation Studies} \label{sec:sim}

\subsection{Testing uniformity} \label{subsec:testingunif}

We first perform simulation studies for the uniformity tests described in Section~\ref{subsec:uniform}. For the sphere $S^2$, we consider three versions of our test:~an asymptotic test based on  Theorem~\ref{th:asy}, and sampling and bootstrap-based tests as described in Section~\ref{subsec:simple}. For the asymptotic test, we consider $\lambda \in \{1/10,1/4,1,4,10\}$ and use the Imhof method to approximate the asymptotic distribution, and in the case of the Monte Carlo tests only $\lambda = 1$ is used. For the Stiefel manifold $V_{3,2}$ we use the two Monte Carlo tests, as described in Section~\ref{subsec:simple}, with $\Lambda = I_2$. Further, all Monte Carlo tests are based on $K = 199$ independent random samples. We also compare our uniformity tests with the Rayleigh test, the projected Rothman test ($PR_t$, with parameter $t = 0.333$), the projected Cram\'er-von Mises test (PCvm), the projected Anderson-Darling test (PAD), and Gin\'{e}'s $F_n$ and $G_n$ tests. All competitors were implemented using the R package \texttt{sphunif} \cite{sphunif-package}.

To investigate the Type I error rate for the nominal level $0.05$, we drew $N$ independent random samples (where $N=10,000$ for $S^2$, and $N=1,000$ for $V_{3,2}$), each of size $n$, from a uniform distribution.  The empirical Type I error rate is then calculated as the fraction of samples for which the p-value does not exceed $0.05$, the results are displayed in Table \ref{tab:type1_unif}. We note that, as expected, the sampling-based test works well even for very small sample sizes. The bootstrap-based test on $S^2$ performs well for moderate sample sizes $n \ge 50$, but is markedly unconservative for $n = 3$ and $n = 5$. The bootstrap-based test on $V_{3,2}$ exhibits conservative behavior for all sample sizes smaller than $500$. Except for $\lambda = 10$, the asymptotic tests controls the nominal Type I error rate in all settings; for $n \ge 20$, the empirical Type I error rate is very close to the nominal level.

To investigate the power of the test, we drew $n_{sim}$ independent random samples (where $n_{sim}=10,000$ for $S^2$, and $n_{sim}=1,000$ for $V_{3,2}$), each of size $n$, from two non-uniform scenarios. These two scenarios are: 
\begin{itemize}
    \item[(S1)] Data are drawn from a matrix Fisher distribution with parameters $A_1= \begin{pmatrix} 0 \\ 3/5 \\4/5\end{pmatrix}$ for the sphere $S^2$, and $A_1= \frac{1}{2}\begin{pmatrix} 0 & 1 \\ 3/5 & 0 \\4/5 & 0\end{pmatrix}$ for the Stiefel manifold $V_{3,2}$.
    \item[(S2)] $50 \%$ of the data are drawn from a matrix Fisher distribution with parameter $A_{21} = c A_1$, and the remaining $50 \%$ of the data are drawn from a matrix Fisher distribution with $A_{21} = -c A_1$. Here, we set $c=3$ for $S^2$ and $c=2 \sqrt{2}$ for $V_{3,2}$.
\end{itemize}
The results for the power are displayed in Table \ref{tab:power_uniform}. For $S^2$, the Rayleigh test outperforms all other procedures in Scenario (S1), whereas Gin\'e's $G_n$ shows the best performance for Scenario (S2). This is unsurprising, since the Rayleigh test is the uniformly most powerful test for Fisher alternatives on the sphere and $G_n$ is tailored for axial (i.e., antipodal) alternatives. Our tests with $\lambda \in \{1/10,1/4\}$ show results very similar to the Rayleigh test, demonstrating the fact that our test converges to the Rayleigh test as $\lambda \to 0$, cf., Remark \ref{rem:limiting_test_statistic}.

Only a few tests show very good performance for both scenarios:~PAD, Gin\'e's $F_n$, and our tests with $\lambda \in \{1,4,10\}$. The best overall performance is achieved by $\lambda = 1$ and $\lambda =4$, which motivates the use of those parameter values for the real-world data examples. For $V_{3,2}$, the sampling-based test markedly outperforms the bootstrap test, particularly for smaller sample sizes.
\begin{table}[ht]
\caption{Empirical Type I error rates for tests of uniformity.  Values are given in bold if the corresponding Clopper-Pearson confidence interval contains the nominal level $0.05$. Values are given in red if the lower bound of the corresponding Clopper-Pearson confidence interval is greater than the nominal level $0.05$.}
\label{tab:type1_unif}
\centering
\begin{tabular}{l|cccccccc}
  & \multicolumn{8}{c}{$n$} \\
  & 3 & 5 & 10 & 20 & 50 & 100 & 200 & 500 \\
  \hline
\textbf{Sphere $S^{2}$} \\
Asymp. ($\lambda=1$) & 0.033 & 0.040 & 0.043 & \textbf{0.046} & \textbf{0.052} & \textbf{0.051} & \textbf{0.047} & \textbf{0.046} \\ 
Sampling ($\lambda=1$) & \textbf{0.048} & \textbf{0.050} & \textbf{0.048} & \textbf{0.050} & \textbf{0.052} & \textbf{0.051} & \textbf{0.049} & \textbf{0.047} \\ 
Bootstrap ($\lambda=1$) & \textcolor{red}{0.169} & \textcolor{red}{0.069} & \textbf{0.049} & 0.044 & \textbf{0.049} & \textbf{0.049} & \textbf{0.048} & \textbf{0.048} \\ 
Asymp. ($\lambda=1/4$) & 0.025 & 0.039 & 0.044 & \textbf{0.047} & \textbf{0.052} & \textbf{0.049} & \textbf{0.049} & \textbf{0.050} \\ 
Asymp. ($\lambda=4$) & \textbf{0.047} & \textbf{0.048} & \textbf{0.047} & \textbf{0.048} & \textbf{0.051} & \textbf{0.050} & \textbf{0.049} & \textbf{0.048} \\ 
Asymp. ($\lambda=1/10$) & 0.021 & 0.038 & 0.043 & \textbf{0.047} & \textbf{0.052} & \textbf{0.049} & \textbf{0.050} & \textbf{0.051} \\ 
Asymp. ($\lambda=10$) & \textbf{0.052} & \textcolor{red}{0.064} & \textbf{0.052} & \textbf{0.051} & \textcolor{red}{0.056} & \textbf{0.052} & \textbf{0.049} & \textbf{0.048} \\ 
Rayleigh & 0.017 & 0.038 & 0.042 & \textbf{0.047} & \textbf{0.052} & \textbf{0.050} & \textbf{0.050} & \textbf{0.051} \\ 
PR$_t$ & 0.025 & 0.039 & \textbf{0.045} & \textbf{0.047} & \textbf{0.051} & \textbf{0.050} & \textbf{0.050} & \textbf{0.050} \\ 
PCvM & 0.027 & 0.038 & 0.044 & \textbf{0.047} & \textbf{0.052} & \textbf{0.049} & \textbf{0.050} & \textbf{0.049} \\ 
PAD & 0.029 & 0.039 & 0.044 & \textbf{0.047} & \textbf{0.053} & \textbf{0.050} & \textbf{0.049} & \textbf{0.048} \\ 
Gin\'e's $G_{n}$ & 0.032 & 0.039 & \textbf{0.046} & \textbf{0.050} & \textbf{0.050} & \textbf{0.051} & \textbf{0.051} & \textbf{0.052} \\ 
Gin\'e's $F_{n}$ & 0.031 & 0.039 & 0.044 & \textbf{0.046} & \textbf{0.051} & \textbf{0.050} & \textbf{0.048} & \textbf{0.047} \\

\\
\hline
\textbf{Stiefel manifold $V_{3,2}$} \\
Sampling ($\Lambda=I_2$) & \textbf{0.044} & \textbf{0.041} & \textbf{0.062} & \textbf{0.039} & \textbf{0.050} & \textbf{0.063} & \textbf{0.047} & \textbf{0.057} \\ 
Bootstrap ($\Lambda=I_2$) & 0.032 & 0.011 & 0.004 & 0.006 & 0.011 & 0.033 & 0.036 & \textbf{0.051} \\ 
\end{tabular}
\end{table}

\begin{table}[ht]
\centering
\caption{Empirical power for tests of uniformity.  Values of the method with the highest power for each scenario are given in bold.}
\label{tab:power_uniform}
\setlength{\tabcolsep}{6pt}
\begin{tabular}{l| l*{4}{c}@{\hspace{1em}}*{4}{c}}
& \multicolumn{4}{c}{Scenario (S1)} && \multicolumn{4}{c}{Scenario (S2)} \\
Method & $n{=}10$ & $n{=}20$ & $n{=}50$ & $n{=}100$ && $n{=}10$ & $n{=}20$ & $n{=}50$ & $n{=}100$ \\
\hline
\multicolumn{10}{l}{\textbf{Sphere $S^{2}$}}\\
Asymp. ($\lambda=1$) & 0.230 & 0.477 & 0.906 & 0.998 && 0.020 & 0.156 & 0.757 & 0.995 \\
Sampling ($\lambda=1$) & \textbf{0.250} & 0.481 & 0.905 & 0.998 && 0.025 & 0.164 & 0.750 & 0.995 \\
Bootstrap ($\lambda=1$) & 0.244 & 0.473 & 0.901 & 0.998 && 0.011 & 0.108 & 0.700 & 0.994 \\
Asymp. ($\lambda=1/10$) & \textbf{0.250} & 0.507 & \textbf{0.922}& \textbf{0.999} && 0.004 & 0.005 & 0.007 & 0.020 \\
Asymp. ($\lambda=1/4$) & \textbf{0.250} & 0.505 & 0.921 & \textbf{0.999} && 0.004 & 0.008 & 0.029 & 0.295 \\
Asymp. ($\lambda=4$) & 0.165 & 0.330 & 0.783 & 0.991 && 0.123 & 0.390 & 0.909 & 0.999 \\
Asymp. ($\lambda=10$) & 0.124 & 0.221 & 0.605 & 0.947 && 0.122 & 0.319 & 0.832 & 0.996 \\
Rayleigh & \textbf{0.250} & \textbf{0.509} & \textbf{0.922} & \textbf{0.999} && 0.003 & 0.004 & 0.003 & 0.003 \\
PR$_t$ & 0.249 & 0.504 & 0.921 & \textbf{0.999} && 0.004 & 0.008 & 0.046 & 0.424 \\
PCvM & 0.248 & 0.502 & 0.920 & \textbf{0.999} && 0.004 & 0.012 & 0.115 & 0.699 \\
PAD & 0.244 & 0.496 & 0.917 & \textbf{0.999} && 0.005 & 0.026 & 0.323 & 0.919 \\
Gin\'e's $F_{n}$ & 0.241 & 0.493 & 0.916 & \textbf{0.999} && 0.007 & 0.051 & 0.507 & 0.973 \\
Gin\'e's $G_{n}$ & 0.049 & 0.058 & 0.091 & 0.149 && \textbf{0.338} & \textbf{0.702} & \textbf{0.989} & \textbf{1.000} \\

\hline
\multicolumn{10}{l}{\textbf{Stiefel manifold $V_{3,2}$}}\\
Sampling ($\Lambda=I_2$) & \textbf{0.098} & \textbf{0.158} & \textbf{0.422} & \textbf{0.792} && \textbf{0.086} & \textbf{0.203} & \textbf{0.660} & \textbf{0.970} \\
Bootstrap ($\Lambda=I_2$) & 0.016 & 0.044 & 0.267 & 0.716 && 0.009 & 0.050 & 0.457 & 0.942 \\
\end{tabular}
\end{table}

\subsection{Testing simple Fisher hypotheses} \label{subsec:simplefisher}

Next we investigate the performance for testing that the distribution follows a matrix Fisher distribution with a fixed parameter $A_0$, where $A_0$ was chosen as
$A_0= \begin{pmatrix} 0 \\ 3/5 \\4/5\end{pmatrix}$ for the sphere $S^2$ and $A_0= \begin{pmatrix} 0 & 1 \\ 3/5 & 0 \\4/5 & 0\end{pmatrix}$ for the Stiefel manifold $V_{3,2}$. For this setting, the sampling and bootstrap-based tests as described in Section \ref{subsec:simple} are used. All Monte Carlo tests are based on drawing $K = 199$ samples; the constant $U_3$ required for the bootstrap test is calculated from a sample of size  $N=50,000$.
The empirical Type I error rate at the nominal level $\alpha =0.05$ is very similar to the empirical Type I error rate of the corresponding Monte Carlo procedures when testing uniformity and is provided in Table 1 of the Supplementary Material.

To investigate the power, we drew $n_{sim}$ independent random samples (where $n_{sim}=10,000$ for $S^2$, and $n_{sim}=1,000$ for $V_{3,2}$), each of size $n$, from the following two distributions: 
\begin{itemize}
    \item[(S1)] Matrix Fisher distribution with parameters $A_1= 1/4 \,\begin{pmatrix} 0 \\ 3/5 \\4/5\end{pmatrix}$ for the sphere $S^2$ and $A_1= \frac{1}{2}\begin{pmatrix} 0 & 1 \\ 3/5 & 0 \\4/5 & 0\end{pmatrix}$ for the Stiefel manifold $V_{3,2}$.
    \item[(S2)] Matrix Fisher distribution with parameters $A_2= \,\begin{pmatrix} 0 \\ 1 \\ 0\end{pmatrix}$ for the sphere $S^2$ and $A_2= \begin{pmatrix} 0 & 1 \\ 1 & 0 \\0 & 0\end{pmatrix}$ for the Stiefel manifold $V_{3,2}$.
\end{itemize}

The first alternative represents a difference in concentration (but not location) compared to the null distribution, whereas the second alternative represents a difference in location (but not in concentration). The results for the empirical power are provided in Table \ref{tab:power_simple}. The results for $\lambda \in \{1,4,10\}$ provide numerical evidence that the test is consistent against both alternatives; as observed before, the sampling-based test outperforms the bootstrap test.

\begin{table}[ht]
\centering
\caption{Empirical power for testing that the distribution follows a Fisher distribution with density $g (\cdot; A_0,0)$.  Values of the method with the highest power for each scenario are given in bold.}
\label{tab:power_simple}
\setlength{\tabcolsep}{6pt}
\begin{tabular}{l| l*{4}{c}@{\hspace{1em}}*{4}{c}}
& \multicolumn{4}{c}{Scenario (S1)} && \multicolumn{4}{c}{Scenario (S2)} \\
Method & $n{=}10$ & $n{=}20$ & $n{=}50$ & $n{=}100$ && $n{=}10$ & $n{=}20$ & $n{=}50$ & $n{=}100$ \\
\hline
\multicolumn{10}{l}{\textbf{Sphere $S^{2}$}}\\
Sampling ($\lambda=1$) & \textbf{0.181} & \textbf{0.306} & \textbf{0.659} & \textbf{0.934} && \textbf{0.219} & \textbf{0.421} & \textbf{0.837} & \textbf{0.992} \\
Bootstrap ($\lambda=1$) & 0.137 & 0.253 & 0.619 & 0.920 && 0.197 & 0.391 & 0.821 & 0.990 \\
\hline
\multicolumn{10}{l}{\textbf{Stiefel manifold $V_{3,2}$}}\\
Sampling ($\Lambda=I_2$) & \textbf{0.095} & \textbf{0.197} & \textbf{0.462} & \textbf{0.808} && \textbf{0.152} & \textbf{0.278} & \textbf{0.686} & \textbf{0.957} \\ 
Bootstrap ($\Lambda=I_2$) & 0.010 & 0.059 & 0.347 & 0.772 && 0.039 & 0.121 & 0.565& 0.942 \\ 
\end{tabular}
\end{table}

\subsection{Testing composite hypotheses} \label{subsec:simcomp}

Finally, we tested the composite hypothesis that a random sample follows a matrix Fisher distribution with an unspecified parameter $A$; in so doing, we applied the parametric bootstrap approach described in Section \ref{subsec:comp}. As estimator $\widehat{A}_n$ for $A$, we used the maximum likelihood (ML) estimator derived in \cite{jupp1979maximum}. The required derivatives of the generalized hypergeometric functions were approximated discretely, and the system of likelihood equations was solved using the package \texttt{nleqslv}. All required $U_3$ terms were approximated using a sample of size $N=1,000$. For the sphere $S^2$, $K=199$ bootstrap samples were used, and for the Stiefel manifold $V_{3,2}$ we used $K=39$ bootstrap samples.
To evaluate the empirical Type I error rate, we drew $n_{sim} = 1,000$ independent random samples of size $n$ from a matrix Fisher distribution with parameter $A_0= \begin{pmatrix} 0 \\ 3/5 \\4/5\end{pmatrix}$ for $S^2$ and $A_0= \begin{pmatrix} 0 & 1 \\ 3/5 & 0 \\4/5 & 0\end{pmatrix}$ for $V_{3,2}$. The results are provided in Table \ref{tab:type1_composite}. The test shows unconservative behavior for very small sample sizes, but controls the Type I error rate satisfyingly for sample sizes of $n \ge 20$. 

To investigate the power, we drew $n_{sim} = 1,000$ independent random samples of size $n$ from the following scenarios:
\begin{itemize}
    \item[(S1)]  $50 \%$ of the data are drawn from a matrix Fisher distribution with parameter $A_{11} = c A_0$, the other $50 \%$ of the data is drawn from a matrix Fisher distribution with $A_{12} = -c A_0$. The constant $c$ is set equal to $3$ for $S^2$ and equal to $2$ for $V_{3,2}$.
    \item[(S2)]  $50 \%$ of the data are drawn from a matrix Fisher distribution with parameter $A_{21} = 6 A_0$, the other $50 \%$ of the data is drawn from a matrix Fisher distribution with $A_{22} = 0.1 A_0$. 
\end{itemize}

\noindent
The results provide numerical evidence that the test is consistent against both alternatives.

\begin{table}[t!]
\caption{Empirical Type I error rates for testing that the distribution follows a Fisher distribution.  Values are given in bold if the corresponding Clopper-Pearson confidence interval contains the nominal level $0.05$. Values are given in red if the lower bound of the corresponding Clopper-Pearson confidence interval is greater than the nominal level $0.05$.}
\label{tab:type1_composite}
\centering
\begin{tabular}{l|cccccccc}
  & \multicolumn{8}{c}{$n$} \\
  & 3 & 5 & 10 & 20 & 50 & 100 & 200 & 500 \\
\hline
\textbf{Sphere $S^{2}$} \\
P. Bootstrap ($\lambda=1$) & \textcolor{red}{0.089} & \textcolor{red}{0.083} & \textcolor{red}{0.076} & \textbf{0.057}  & \textbf{0.049} & \textbf{0.056} & \textbf{0.046} & \textbf{0.042} \\ 
\hline
\textbf{Stiefel manifold $V_{3,2}$} \\
P. Bootstrap ($\Lambda=I_2$) &  \textcolor{red}{0.152} &  \textcolor{red}{0.110} & \textcolor{red}{0.089} & \textbf{0.060} & \textbf{0.062} & \textbf{0.059} & \textbf{0.060} & \textbf{0.052} \\  
\end{tabular}
\end{table}

\medskip

\begin{table}[t]
\centering
\caption{Empirical power for testing that the distribution follows a Fisher distribution.}
\label{tab:power_composite}
\setlength{\tabcolsep}{6pt}
\begin{tabular}{l| l*{4}{c}@{\hspace{1em}}*{4}{c}}
& \multicolumn{4}{c}{Scenario (S1)} && \multicolumn{4}{c}{Scenario (S2)} \\
Method & $n{=}10$ & $n{=}20$ & $n{=}50$ & $n{=}100$ && $n{=}10$ & $n{=}20$ & $n{=}50$ & $n{=}100$ \\
\hline
\multicolumn{10}{l}{\textbf{Sphere $S^{2}$}}\\
P. Bootstrap ($\lambda=1$) & 0.437 & 0.722 & 0.984 & 1.000 && 0.190 & 0.331 & 0.701 & 0.956 \\
\hline
\multicolumn{10}{l}{\textbf{Stiefel manifold $V_{3,2}$}}\\
P. Bootstrap ($\Lambda=I_2$) & 0.626 & 0.881 & 0.998 & 1.000 && 0.468 & 0.745 & 0.989 & 1.000 \\ 
\end{tabular}
\end{table}

\section{An application to data on the orbit of comets} \label{sec:real}

As an application to real-world data we consider the well-known data set, on the orbit of comets, from the Jet Propulsion Laboratory's Small-Body Database Search Engine available at \url{https://ssd.jpl.nasa.gov/tools/sbdb_query.html} and in the R package \texttt{sphunif} on CRAN. To be consistent with previous analyses, we remove the duplicate entries in rows 13-15, yielding the same data set of $n=208$ comets analyzed in \cite{cuesta2009projection}. As in \cite{cuesta2009projection}, we calculate the normal vector of each comet orbit,
\[
\mathbf{n} =
\begin{pmatrix}
\sin i \, \sin \Omega \\
-\sin i \, \cos \Omega \\
\cos i
\end{pmatrix},
\]
where $i$ is the inclination of the orbit.

\begin{figure}
	\centering
	\includegraphics[width=0.7 \textwidth]{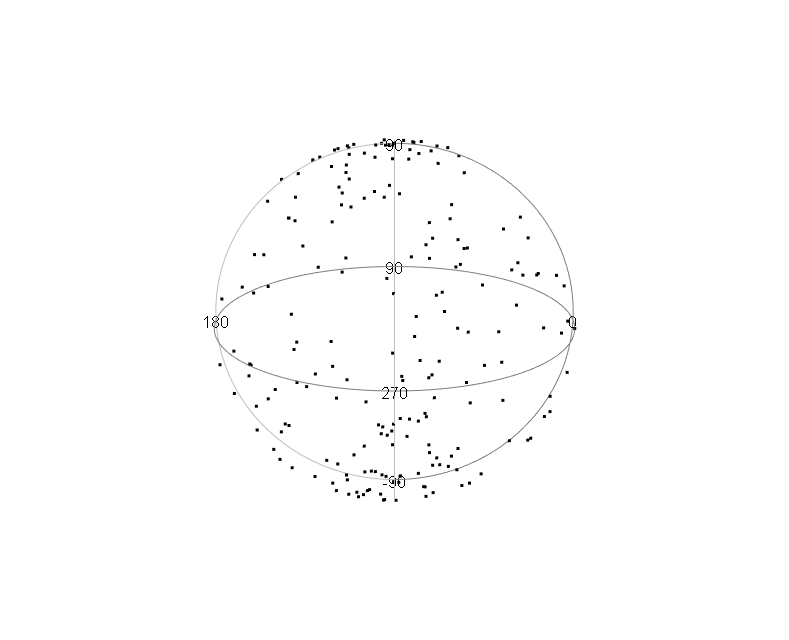}
	\caption{Illustration of the normal vectors for the comet orbits. The reference plane represents the ecliptic.}
	\label{fig:sphere}
\end{figure}

The normal vectors are visualized in Figure \ref{fig:sphere}, the reference plane (``equator'') representing the ecliptic. The figure suggests that there may be an accumulation of orbits near the ecliptic, which is consistent with the analysis in \cite{jupp2003intrinsic}, where it is suggested that such an accumulation likely is due to an observational bias.

We first apply several tests of uniformity to assess whether the normal vectors are uniformly distributed on $S^2$; the results are displayed in Table 6. 
 Many of the available test statistics fail to reject the null hypothesis of uniformity at a nominal level of significance of $0.05$, a notable exception being Gin\'e's statistic, $G_n$, that is explicitly tailored for antipodal distributions.  By contrast, our testing procedures show a significant result despite there being no specialized tailoring toward any alternative distribution.

We also use the parametric bootstrap testing procedure described in Section \ref{subsec:simcomp} to test the null hypothesis that the normal vectors follow a Fisher distribution; both the test with $\lambda = 1$ and the test with $\lambda = 4$ reject the null hypothesis at a level of $\alpha = 0.05$.

\begin{table}[ht]
\centering
\begin{minipage}[t]{0.48\textwidth}
\centering
\caption{P-values of different methods for testing uniformity for the normal vectors of the comet orbits.}
\vspace{5pt}
\begin{tabular}{lc}
\hline
\textbf{Method} & \textbf{P-value}  \\ 
\hline
Asymp. ($\lambda=1$) &  $\mathbf{0.0381}$ \\ 
Bootstrap ($\lambda=1$) &  $\mathbf{0.0371}$\\ 
Sampling ($\lambda=1$)  &  $\mathbf{0.0335}$ \\ 
Asymp. ($\lambda=4$) &  $\mathbf{0.0044}$ \\ 
Bootstrap ($\lambda=4$) &  $\mathbf{0.0075}$ \\ 
Sampling ($\lambda=4$)  &  $\mathbf{0.0050}$ \\ 
Gin\'e's $F_{n}$ & $0.0592$ \\ 
Gin\'e's $G_{n}$ & $\mathbf{0.0086}$   \\ 
PAD &  $0.0693$ \\ 
PCvM &  $0.0947$ \\ 
PR$_t$ &  $0.1116$ \\ 
Rayleigh & $0.2009$   \\ 
\hline
\end{tabular}

\end{minipage}
\hfill
\begin{minipage}[t]{0.48\textwidth}
\centering
\caption{P-values of the proposed parametric bootstrap approach for testing whether the normal vectors of the comet orbits follow a Fisher distribution.}
\vspace{5pt}
\begin{tabular}{lc}
\hline
\textbf{Method} & \textbf{P-value}  \\ 
\hline
P. Bootstrap ($\lambda=1$) &  $\mathbf{0.0044}$\\ 
P. Bootstrap ($\lambda=4$) &  $\mathbf{0.0026}$ \\ 
\hline
\end{tabular}
\label{tab:comet2}
\end{minipage}
\end{table}

Additionally, we calculate the unit vector pointing to the perihelion
\[
\mathbf{p}=
\begin{pmatrix}
\cos \Omega \, \cos \omega - \sin \Omega \, \cos i \, \sin \omega \\[4pt]
\sin \Omega \, \cos \omega + \cos \Omega \, \cos i \, \sin \omega \\[4pt]
\sin i \, \sin \omega
\end{pmatrix},
\]
where $\Omega$ is the longitude of the ascending node and $\omega$ is the argument of the perihelion.

Then $(\mathbf{n},\mathbf{p})$ lies in the Stiefel manifold $V_{3,2}$. We then apply our testing procedures to test the null hypotheses that a) $(\mathbf{n},\mathbf{p})$ follows a uniform distribution on $V_{3,2}$ and b) $(\mathbf{n},\mathbf{p})$ follows a matrix Fisher distribution on $V_{3,2}$. The results are displayed in Table 8 and Table 9, respectively.

\begin{table}[ht]
\centering
\begin{minipage}[t]{0.48\textwidth}
\centering
\label{tab:cometstiefel1}
\caption{P-values of different methods for testing uniformity of $(\mathbf{n},\mathbf{p})$ in $V_{3,2}$, where $\mathbf{n}$ denotes the normal vector of the comet orbit and $\mathbf{p}$ is the vector of the perihelion. }
\vspace{5pt}
\begin{tabular}{lc}
\hline
\textbf{Method} & \textbf{P-value}  \\ 
\hline 
Bootstrap ($\Lambda=I_2$) &  $\mathbf{0.0003}$\\ 
Sampling ($\Lambda=I_2$)  &  $\mathbf{0.0001}$\\ 
Bootstrap ($\Lambda=4 \, I_2$) &  $\mathbf{0.0044}$ \\ 
Sampling ($\Lambda=4 \, I_2$)  &  $\mathbf{0.0001}$\\ 
\hline
\end{tabular}
\label{tab:comet1}
\end{minipage}
\hfill
\begin{minipage}[t]{0.48\textwidth}
\centering
\label{tab:cometstiefel2}
\caption{P-values of the proposed parametric bootstrap approach for testing whether $(\mathbf{n},\mathbf{p})$ follow a matrix Fisher distribution in $V_{3,2}$.}
\vspace{5pt}
\begin{tabular}{lc}
\hline
\textbf{Method} & \textbf{P-value}  \\ 
\hline
P. Bootstrap ($\Lambda=I_2$) &  $\mathbf{0.0009}$\\ 
P. Bootstrap ($\Lambda=4 \, I_2$) &  $\mathbf{0.0005}$ \\ 
\hline
\end{tabular}
\label{tab:second}
\end{minipage}
\end{table}

\section{Discussion and outlook} \label{sec:discussion}

In this article, we have derived a general framework for goodness-of-fit testing on Stiefel manifolds. The proposed tests are consistent against all alternatives and demonstrate very good performance across a broad range of scenarios in our simulation studies.

The framework is connected naturally to the classical distance- and kernel-based testing literature (see, e.g., \cite{rizzo2016energy, sejdinovic2013equivalence, edelmann2022regression}). This connection also indicates generalizations of the results given here to two-sample testing and independence testing on Stiefel manifolds, which are topics that we intend to develop in future work.

In other directions, the results developed in the present paper may also be extended to distributions beyond the Fisher-Bingham family. Exploring such generalizations will enable applications of our testing framework to a wider range of statistical models arising in manifold-valued data analysis.

Finally, a promising but substantially more challenging direction concerns the development of a complete extension of the Funk-Hecke theorem to Stiefel manifolds; here, we remark that partial extensions of that theorem have been treated by several authors, e.g., \cite{rubin2024injectivity}. The full extension of the Funk-Hecke theorem would enable us to derive the complete asymptotic distribution of our test statistics in the general Stiefel manifold setting, thereby extending Theorem \ref{th:asy}.


\end{document}